\documentclass[lettersize,journal]{IEEEtran}
\usepackage{amsmath,amsfonts}
\usepackage{amsthm}
\newtheorem{proposition}{Proposition}
\usepackage{algorithm}
\usepackage{array}
\usepackage[caption=false,font=normalsize,labelfont=sf,textfont=sf]{subfig}
\usepackage{textcomp}
\usepackage{stfloats}
\usepackage{url}
\usepackage{verbatim}
\usepackage{graphicx}
\usepackage{color}
\graphicspath{\figures}
\usepackage{wrapfig}
\usepackage{cite}
\usepackage{multirow}

\usepackage{algpseudocode}

\begin{document}

\title{A Scalable Network-Aware Multi-Agent Reinforcement Learning Framework for Decentralized Inverter-based Voltage Control}

\author{Han Xu,~\IEEEmembership{Student Member,~IEEE,} Jialin Zheng,~\IEEEmembership{Student Member,~IEEE,} Guannan Qu,~\IEEEmembership{Member,~IEEE}
}



\maketitle

\begin{abstract}
 This paper addresses the challenges associated with decentralized voltage control in power grids due to an increase in distributed generations (DGs). Traditional model-based voltage control methods struggle with the rapid energy fluctuations and uncertainties of these DGs. While multi-agent reinforcement learning (MARL) has shown potential for decentralized secondary control, scalability issues arise when dealing with a large number of DGs. This problem lies in the dominant centralized training and decentralized execution (CTDE) framework, where the critics take global observations and actions. To overcome these challenges, we propose a scalable network-aware (SNA) framework that leverages network structure to truncate the input to the critic's Q-function, thereby improving scalability and reducing communication costs during training. Further, the SNA framework is theoretically grounded with provable approximation guarantee, and it can seamlessly integrate with multiple multi-agent actor-critic algorithms. The proposed SNA framework is successfully demonstrated in a system with 114 DGs, providing a promising solution for decentralized voltage control in increasingly complex power grid systems.
\end{abstract}


\section{Introduction}
\IEEEPARstart{R}{newable} energy plays a vital role in mitigating climate change by reducing greenhouse gas emissions \cite{osmanCostEnvironmentalImpact2023}. The rising utilization of renewable energy has catalyzed the transformation of power grids, leading to an increasing number of distributed generations (DGs). As a result, the rapid energy fluctuations and uncertainties of these DGs present challenges to the voltage control problem \cite{sunReviewChallengesResearch2019} where conventional voltage regulation methods based on capacitors and voltage regulators \cite{hammadComparingVoltageControl1996,attarNovelStrategyOptimal2018} may fail to maintain the voltage magnitude within an acceptable range. Fortunately, given that a majority of DGs are inverter-based DGs, they possess the capability to deliver rapid voltage regulation utilizing their control flexibility \cite{lasseterGridFormingInvertersCritical2020}. 

For controlling DGs, there are two levels of control including primary and secondary control \cite{bidramHierarchicalStructureMicrogrids2012}. The primary control usually adopts droop technique to maintain system stability and mitigate the voltage fluctuation in a localized manner \cite{sunNewPerspectivesDroop2017}. However, solely using primary control can cause voltage deviations \cite{bidramSecondaryControlMicrogrids2013}. To maintain the voltage at the predefined nominal value, secondary control is necessary \cite{olivaresTrendsMicrogridControl2014,rajeshReviewControlAc2017}. 
In this paper, we focus on \emph{decentralized} secondary control. Compared with centralized secondary control, decentralized secondary control holds advantages, including lower communication burdens and easier deployments \cite{vovosCentralizedDistributedVoltage2007}, as in decentralized settings, these agents make decisions that collaboratively mitigate the system voltage fluctuation with only \emph{local measurements}. 
The design of decentralized control in general has long been a challenge in control theory \cite{blondelSurveyComputationalComplexity2000a}.  
Further complicating the controller design is the large number of DGs, the nonlinear dynamics of power grids and the uncertainty associated with loads and renewable energies. 
All these make the design of model-based decentralized secondary control challenging, thereby negatively impacting the control performance \cite{dehkordiDistributedNoiseResilientSecondary2019,simpson-porcoSecondaryFrequencyVoltage2015}.

Reinforcement learning (RL), a model-free approach, has been recognized to have great potential for realizing secondary control and other power system control problems due to its adaptability and ability to solve complex control problems \cite{chenMultiAgentReinforcementLearning2022}. Our decentralized secondary control falls under the realm of multi-agent RL (MARL). In MARL, in order to make these different DGs (e.g., agents) collaborate together to maintain the system voltage, a popular framework called Centralized Training and Decentralized Execution (CTDE) is typically adopted \cite{caneseMultiAgentReinforcementLearning2021,chenPowerNetMultiAgentDeep2022c,yangDistributedDynamicInertiadroop2022,liuOnlineMultiAgentReinforcement2021}. In this framework, the critics take global information and measure how individual actions influence global welfare, and therefore guide actors to update policies to increase global welfare. 

Unfortunately, this dominant framework suffers from poor scalability. As mentioned above, in this framework, the Q-function of the critics needs to take global observations and actions to estimate expected rewards and train individual actors. Therefore, the size of the observation and action spaces of the critics exponentially increases with the number of agents, rendering the lack of scalability of this framework \cite{caneseMultiAgentReinforcementLearning2021,quScalableReinforcementLearning2021}. More specifically, when the number of DGs (e.g., agents) is large and the DGs are controlled in a decentralized manner, this framework suffers from low sample efficiency and may fail to train effective critics and actors to provide desirable control performance due to the notorious curse of dimensionality. For this reason, the number of agents in the existing literature on decentralized voltage control can only reach up to 40 \cite{chenPowerNetMultiAgentDeep2022c}, which is far below realistic networks. To address this, some methods try to utilize advanced neural networks to improve the scalability. In \cite{yeScalablePrivacyPreservingMultiAgent2021}, the multi-actor-attention-critic (MAAC) algorithm is used, which utilizes the advanced multi-head attention neural network \cite{iqbalActorAttentionCriticMultiAgentReinforcement}. In \cite{zhouScalableMultiregionPerimeter2023,wangRealTimeJointRegulations2023}, the QMIX method is used. However, these methods still face limitations, as attention networks face the problem of high computational demand due to the square relation between its attention computation and the length of the sequence. Further, QMIX is a type of Q-learning method that only works for finite-action space. In addition, all these methods lack provable guarantees in their ability to handle large systems.

In summary, the conventional CTDE framework encounters scalability problems making it challenging to effectively handle systems with a large number of agents. A common theme in these methods is that they do not leverage the network structure available in this voltage control problem. In contrast, the voltage control problem has rich network structure properties which describe the interaction patterns among agents (e.g., DGs). This network structure can be a valuable source of prior knowledge that can help improve the scalability.

In this paper, we propose a Scalable Network-Aware (SNA) framework that leverages network structure to achieve scalability. Critically, our framework truncates the input to the critic's Q-function, limiting the input to only local and neighboring information, thereby avoiding the need for global information and improving scalability. This means that even for systems with a large number of DGs, high sample efficiency and good training outcomes can be assured. Furthermore, because training only involves information exchange between neighbors, the training can be conducted in a distributed manner, which reduces communication costs. Lastly, our framework is theoretically grounded as we prove the exponential decay property, justifying the truncation.

The main contributions of this paper are summarized as follows:
\begin{enumerate}
  \item The voltage control problem is modeled as a networked Markov Decision Process (MDP), and the reasonableness of this approach is justified. The corresponding training environment has been developed, which can directly read standard matpower data to achieve dynamic simulation.
  \item The proposed SNA framework improves the scalability which is theoretically guaranteed and experimentally validated. Further, our framework is versatile and can seamlessly integrate with different actor-critic implementations, including multi-agent soft-actor-critic (MASAC), multi-agent twin delayed deep deterministic policy gradient (MATD3).
  \item With the proposed framework, voltage control for a system containing 114 DGs (e.g., agents) is successfully realized. To the best of our knowledge, this number exceeds the number of existing literature \cite{chenPowerNetMultiAgentDeep2022c}.
\end{enumerate}

The remainder of this paper is organized as follows. Section II introduces the voltage control problem and its formulation in networked MDP. Section III introduces preliminaries of MARL and the CTDE framework. Section IV presents the proposed SNA framework, combines it with the MASAC algorithm and justifies it with theoretical analysis. Section V demonstrates case studies and result analyses, and conclusions are drawn in Section VI. 

\section{Problem Formulation}

\begin{figure*}[!t]
\centering
\includegraphics[width=7in]{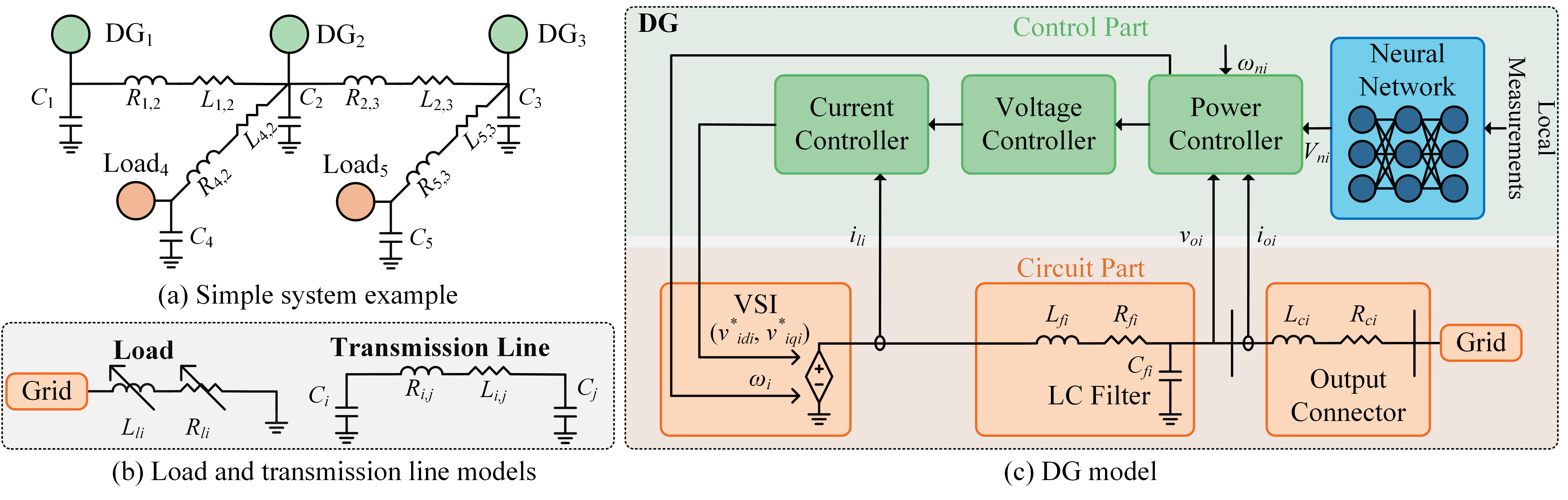}
\caption{Diagram of DGs, loads and transmission lines.}
\label{DER_load_transmission_model}
\end{figure*}

In this section, the voltage control problem is introduced first. Then this problem is formulated as a networked MDP with the justification of why this problem fits in this formulation.

\subsection{Inverter-based Voltage Control Problem}

From the perspective of modelling, there are three basic components in power grids including loads, inverter-based DGs, and transmission lines. Loads and DGs are connected through transmission lines, as depicted in Fig. \ref{DER_load_transmission_model} (a).  The models of loads and transmission lines are shown in Fig. \ref{DER_load_transmission_model} (b). The power consumed by loads typically has stochastic variations and are modelled as time-varying impedance. Apart form the typical loads, non-dispatchable DGs that generate power based on weather conditions are also modeled as negative stochastic loads \cite{houPriorityDrivenSelfOptimizingPower2022}. 

In terms of the model of inverter-based DGs, they consist of the circuit and control part, as shown in Fig. \ref{DER_load_transmission_model} (c). The details of this model can be found in \cite{bidramDistributedCooperativeSecondary2013a}. In the circuit part, the voltage source inverter (VSI) is connected to the remaining grid through an LC filter and output connector. In the control part, the primary control is usually used to ensure voltage and frequency stability. Primary control usually adopts classic droop control strategy which can realize distributed control but alone cannot restore the voltage to predefined nominal value \cite{bidramDistributedCooperativeSecondary2013a,sunNewPerspectivesDroop2017}. Therefore, it necessitates the use of secondary control to compensate for voltage and frequency deviations caused by primary control and load variations. The secondary control is responsible for establishing the voltage reference for the primary control, aiming to regulate both the frequency and voltage amplitude to their predefined nominal values \cite{guerreroHierarchicalControlDroopControlled2011,bidramSecondaryControlMicrogrids2013}. 

In our framework, a neural network-based controller is used to realize the secondary control. It uses only the local measurements and feeds voltage references $V_{mi}$ to the primary controller to maintain the voltage magnitude of the system, as shown in Fig. \ref{DER_load_transmission_model} (c). In summary, the ultimate objective is to design decentralized secondary controllers to maintain the voltage of all inverter-based DGs with the presence of stochastic load variations and imperfections of primary control.

\subsection{Formulation in Networked MDP}

In this section, the above voltage control problem is formulated as a Networked MDP as in \cite{quScalableReinforcementLearning2021}. We start by introducing the general networked MDP, and then show why the voltage control problem fits into this formulation. 

We consider a network of $N$ agents that are associated with an underlying undirected graph $\mathcal{G}=(\mathcal{N},\mathcal{E})$, where $\mathcal{N}=\{1,\dots N\}$ is the set of agents and $\mathcal{E}\subset \mathcal{N}\times \mathcal{N}$ is the set of edges. Each agent $i$ is associated with state $\mathbf{s}_i\in \mathcal{S}_i$, $\mathbf{a}_i\in \mathcal{A}_i$. The global state is denoted as $\mathbf{s} = (\mathbf{s}_1, \dots , \mathbf{s}_N) \in \mathcal{S} := \mathcal{S}_1 \times \dots \times \mathcal{S}_N$ and similarly the global action $\mathbf{a}=(\mathbf{a}_1,\dots, \mathbf{a}_N)\in \mathcal{A}:=\mathcal{A}_1 \times \dots \times \mathcal{A}_N$. At time $t$, given current state $\mathbf{s}_t$ and action $\mathbf{a}_t$, the next individual state $\mathbf{s}_{t+1}$ is independently generated and is only dependent on neighbors:
\begin{equation}
\label{Transition Equation}
P(\mathbf{s}_{t+1}|\mathbf{s}_t,\mathbf{a}_t)=\prod_{i=1}^NP(\mathbf{s}_{i,t+1}|\mathbf{s}_{\mathcal{N}_i,t},\mathbf{a}_{i,t}),
\end{equation}
where notation $\mathcal{N}_i$ means the neighborhood of $i$ (including $i$ itself) and notation $\mathbf{s}_{\mathcal{N}_i,t}$ means the states of the agents in $\mathcal{N}_i$. In addition, for integer $\kappa \ge 0$, we use $\mathcal{N}_i^{\kappa}$ to denote the $\kappa$-hop neighborhood of $i$, i.e. the nodes whose graph distance to i has length less than or equal to $\kappa$, and $\mathcal{N}_{-i}^\kappa=\mathcal{N}/\mathcal{N}_i^\kappa$.

In each time step $t$, each agent $i$ can only observe the local observation $\mathbf{o}_{i,t} \in \mathcal{O}_i \subset \mathcal{S}_i$ rather than the local states. Then each agent chooses its own action $\mathbf{a}_{i,t}$  using a stochastic policy defined as a probability density function $\pi_i : \mathcal{O}_i \times \mathcal{A}_i \rightarrow [0,\infty)$, i.e., $\mathbf{a}_{i,t}\sim \pi_i(\cdot|\mathbf{o}_{i,t})$. 

Furthermore, each agent is associated with a stage reward function $r_{i,t} = r_i(\mathbf{o}_{i,t},\mathbf{a}_{i,t})$ that depends on the local observations and actions, and the global stage reward is $r_t = r(\mathbf{o}_t,\mathbf{a}_t)=\frac{1}{N}\sum_{i=1}^N r_{i,t}$. The objective is to find local policies $\pi_i(\cdot|\mathbf{o}_{i,t})$ such that the discounted global stage reward with discount factor $\gamma$ is maximized.
\begin{equation}
\underset{\pi_1,\cdots,\pi_N}{\text{max}}J(\theta):=\mathbb{E}[\sum_{t=0}^{\infty}\gamma^tr_t].
\end{equation}

In this voltage control problem, each DG is considered as an individual agent. Without loss of generality, a node connected with only loads can also be considered as an agent, albeit lacking the capacity for control. The definitions of state, action, and reward of this specific problem are given below.

\subsubsection{State space} The individual states of an agent are defined as all the state variables of a DG.

\subsubsection{Observation space} The individual observations of an agent are the local measurements. They are defined as $(P_i,Q_i,\delta_i,i_{ld,i},i_{lq,i},i_{od,i},i_{oq,i},v_{od,i},v_{oq,i})$, where $P_i$, $Q_i$ are the measured active and reactive power; $\delta_i$ is the voltage angle referred to the common reference angle; $i_{ld,i},i_{lq,i},i_{od,i},i_{oq,i},v_{od,i},v_{oq,i}$ are current and voltage values in the $dp$ reference frame that are measured from the circuit part of a DG, as labeled in Fig. \ref{DER_load_transmission_model}(c).

\subsubsection{Action space} The individual action is the voltage reference $V_{mi}$ fed to the local primary controller. 

\subsubsection{Reward function} As mentioned above, the individual reward depends only on the local observations and actions. Since the objective is to maintain the voltage, the individual reward function is defined as follows 
\begin{equation}
\label{Individual reward function}
r_{i,t} = 
\begin{cases}
    0.05-|1-v_{o,i}| & \text{if } |1-v_{o,i}|\le 0.05\\
    -|1-v_{o,i}| & \text{if } 0.05 \le |1-v_{o,i}|\le 0.2 \\
    -10 & \text{otherwise}
\end{cases}
,
\end{equation}
where $v_{o,i} = \sqrt{v_{od,i}^2+v_{oq,i}^2}$ is the bus voltage magnitude.

\subsection{Justification for Networked MDP Formulation}

\begin{figure}[!t]
\centering
\includegraphics[width=3.0in]{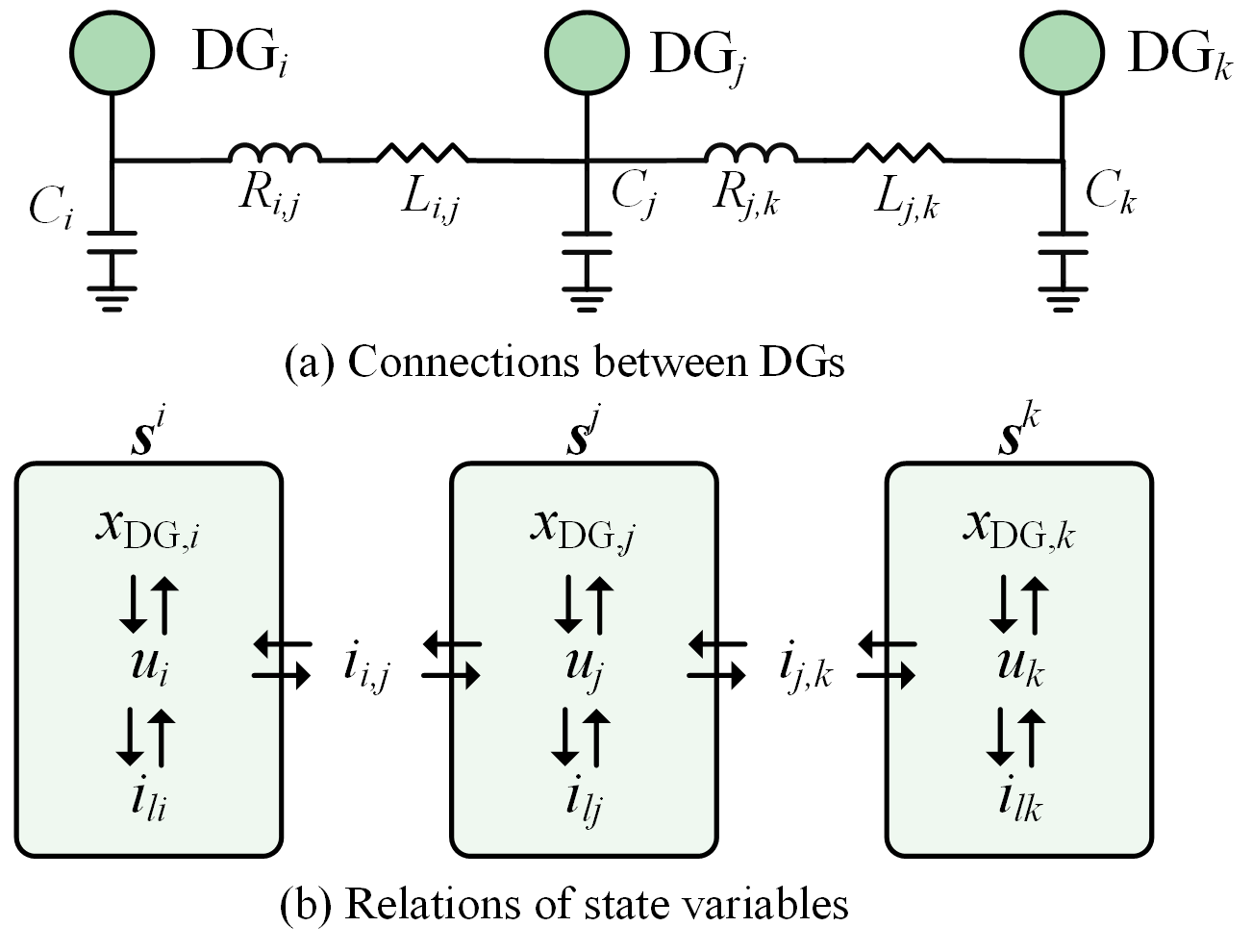}
\caption{Relations between DGs' state variables.}
\label{derivatives_relations}
\end{figure}

In this section, we explain why the voltage control problem can fit in the aforementioned networked MDP. Intuitively, the transition function (\ref{Transition Equation}) is correct due to the delay effect of transmission lines \cite{christopoulosTransmissionLineModelingTLM2022}. This intuition can be further substantiated by examining the relations between the state variables within the system. 

As shown in Fig. \ref{derivatives_relations}(a), we consider a simple system with three DGs. Fig. \ref{derivatives_relations}(b) depicts the interactions among the state variables in the system where $u_i$, $u_j$, $u_k$ are capacitor voltages, $i_{li}$, $i_{lj}$, $i_{lk}$ are load currents, and $x_{\text{DG},i}$, $x_{\text{DG},j}$, $x_{\text{DG},k}$ are state variables of the DGs. The directed arrows in the figure represent that the value of the state variable corresponding to the starting point of the arrow influences the first order derivative of the state variable at the arrow's end point. 

When we examine the state variables of DG$_i$ and DG$_k$ in Fig. \ref{derivatives_relations}, the minimum graph distance between them is four. This implies that a change in the value of any state variables in one of the DGs will only affect derivatives of the fourth order or higher in the other DG. Consequently, if the system is discretized using any numerical method with a precision less than the fourth order, it would satisfy the assumption of the networked MDP, i.e., the interaction between nodes occurs in a local manner. Given the fact that the popular numerical integration method used for electromagnetic transient simulation is trapezoidal method \cite{longEMTPaPowerfulTool1990}, a second-order method, we can conclude the networked MDP modelling is accurate enough from the perspective of numerical modeling and computation.

\section{Preliminaries}

In this section, we introduce preliminaries of MARL and the CTDE framework. Then, we explain the non-scalability problem of the CTDE framework.

\subsection{Multi-agent Reinforcement Learning}

RL, as a model-free approach, has been recognized to have great potential for controlling complex systems. More specifically, for voltage control tasks involving multiple DGs as in our setting, MARL is needed. Deep MARL algorithms based on A2C\cite{chenPowerNetMultiAgentDeep2022c}, DDPG\cite{hossainVoltVAROptimizationDistribution2022}, TD3\cite{wangRealTimeJointRegulations2023}, SAC\cite{yangDistributedDynamicInertiadroop2022} have been adopted in the applications for power system control. These algorithms all belong to multi-agent actor-critic algorithms. Multi-agent actor-critic is an extension of the actor-critic algorithms for environments where multiple agents interact. In such settings, each agent has its own actor and critic. All these actors and critics are approximated by deep neural networks. The actors are responsible for selecting actions, and critics are responsible for evaluating the actions taken by the actors by learning the Q-functions. Through training these actors and critics alternatively, the well-trained critics can gradually help actors improve their policies. The majority of these multi-agent actor-critic algorithms is trained via the CTDE framework, which is discussed in the next section.

\subsection{CTDE Framework}
The prevalent approach in training the above MARL algorithms is the CTDE framework. Within this framework, each agent $i$ adopts a distributed policy $\pi(\mathbf{a}_{i,t}|\mathbf{o}_{i,t})$ to choose an action $\mathbf{a}_{i,t}$ based on its own local observation $\mathbf{o}_{i,t}$. The critics, however, collectively utilize the global observations $\mathbf{o}_t$ and global actions $\mathbf{a}_t$, computing $Q_i(\mathbf{o}_t,\mathbf{a}_t)$ to estimate the global rewards and evaluate how each agent's action influence the global welfare. As a result, agents can learn to cooperatively mitigate voltage violations. Depending on which MARL algorithm is adopted, this framework contains variants like MASAC, MADDPG, MATD3.

In the CTDE framework, critics take the global observations and actions. With a growing number of agents, the notorious ``curse of dimensionality'' tends to manifest. From a theoretical angle, even when the state/action space is finite, the size of the state space is exponentially large in the number of agents which is clearly unscalable. In practice, the Q-functions are approximated by neural networks, in which case the high dimensionality leads to low sample efficiency and makes it difficult to train effective critics as the difficulty of extracting effective features to represent the correct Q-values \cite{otaCanIncreasingInput2020}.
Eventually, the performance of the trained actors may not be satisfactory due to the lack of effective critics. 

\section{Methods}
In this section, we propose the SNA multi-agent deep reinforcement learning framework. In the previous section, the voltage control problem is formulated as a networked MDP which preserves the network structure information. Subsequently, this framework leverages the preserved network structure to achieve superior scalability. 
The main concept is to utilize the characteristic that interactions between agents within a networked MDP occur locally. We demonstrate that this local interaction characteristic implies the exponential decay property, which allows the critics to achieve scalability by performing truncation on the input based on the network structure. Importantly, the truncation will only result in a very small approximation error because of the exponential decay property. 
Then, we combine the proposed SNA framework with the MASAC algorithm and present the pseudo-code. Finally, we provide a theoretical justification for the proposed framework.

\subsection{SNA Framework}

\begin{figure*}[t]
\centering
\includegraphics[width=5.5in]{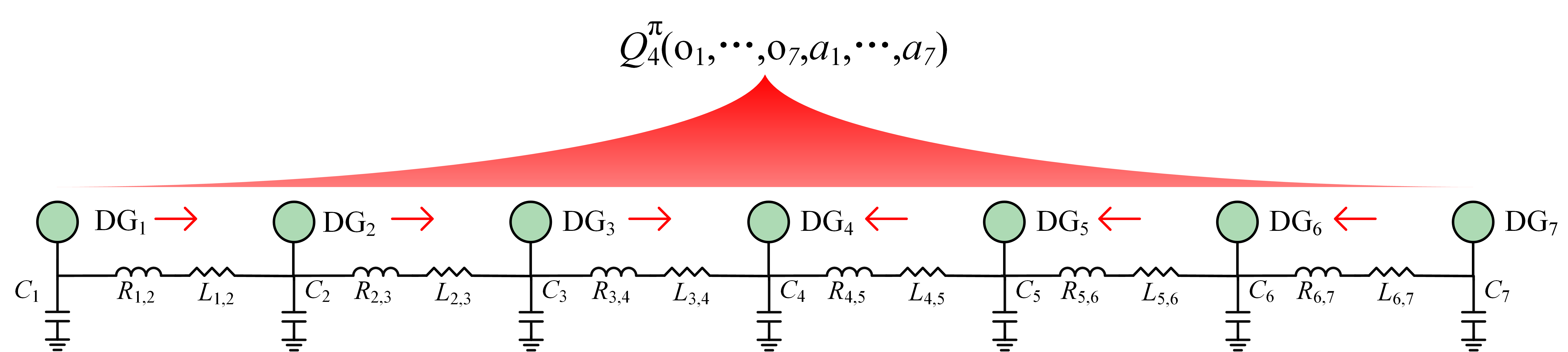}
\caption{Exponential decay property.}
\label{Exponential Decay Property}
\end{figure*}

\begin{figure}[!t]
\centering
\includegraphics[width=3.25in]{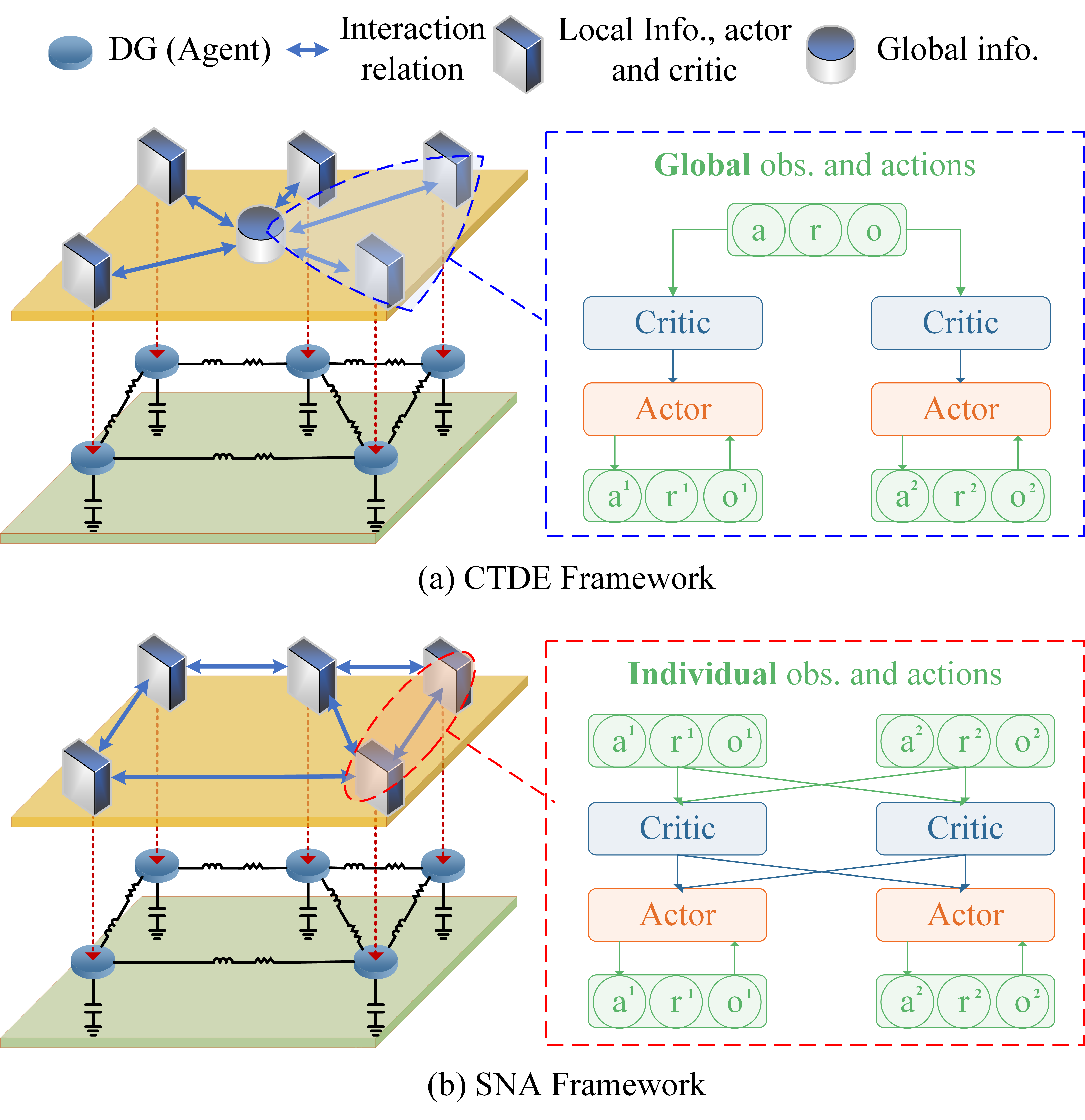}
   \caption{Comparison between CTDE and SNA framework.}
\label{CTDE framework}
\end{figure}

\subsubsection{Exponential decay property}

Here we introduce the exponential decay property. Firstly, the definitions of global and individual Q-functions are given below, and their properties that are related to the network structure is discussed
\begin{equation}
\label{Definition of Q-functions}
\begin{split}
Q^{\pi}(\mathbf{o},\mathbf{a}) = &\mathbb{E}[\sum^{\infty}_{t=0}\gamma^tr_t|\mathbf{o}_0 = \mathbf{o},\mathbf{a}_0 = \mathbf{a}]\\
=&\frac{1}{N}\sum_{i=1}^N\mathbb{E}[\sum_{t=0}^{\infty}\gamma^tr_t^i|\mathbf{o}_0 = \mathbf{o},\mathbf{a}_0=\mathbf{a}]\\
:=&\frac{1}{N}\sum_{i=1}^NQ_i^{\pi}(\mathbf{o},\mathbf{a}),
\end{split}
\end{equation}
where $Q^{\pi}$ is the global Q-function that evaluates the global expected rewards, and $Q_i^{\pi}$ are individual Q-functions that evaluate individual expected rewards. 

As shown in (\ref{Definition of Q-functions}), the global and individual Q-functions are both functions of global observations and actions. Notice that in a networked MDP, interactions between agents happen in a localized manner. 
Intuitively, for a specific agent $i$, given the network structure, it is natural to believe the impact of other agents' observations and actions on agent $i$ decays with increasing topological distance. To formally study this intuition, we define the ''exponential decay property'' of individual Q-functions. The $(c,\rho)$-\textbf{exponential decay property} is said to hold if, for any localized policy $\pi_i$, for any  $i\in \mathcal{N}$, $\mathbf{o}_{\mathcal{N}_i^\kappa}\in \mathcal{O}_{\mathcal{N}_i^\kappa}$, 
$\widetilde{\mathbf{o}}_{\mathcal{N}_{-i}^\kappa}\in \mathbf{O}_{\mathcal{N}_{-i}^\kappa}$, 
$\mathbf{a}_{\mathcal{N}_i^\kappa}\in \mathcal{A}_{\mathcal{N}_i^\kappa}$, $\widetilde{\mathbf{a}}_{\mathcal{N}_{-i}^\kappa}\in \mathcal{A}_{\mathcal{N}_{-i}^\kappa}$, 
$Q_i^{\pi}$ satisfies
\begin{equation}
\label{Eq: Exponential decay property}
\begin{split}
|Q_{i}^{\pi}
(\mathbf{o}_{\mathcal{N}_{i}^\kappa},
\mathbf{o}_{\mathcal{N}_{-i}^\kappa},
\mathbf{a}_{\mathcal{N}_{i}^\kappa},
\mathbf{a}_{\mathcal{N}_{-i}^\kappa}
)-
Q_{i}^{\pi}
(\mathbf{o}_{\mathcal{N}_{i}^\kappa},
\widetilde{\mathbf{o}}_{\mathcal{N}_{-i}^\kappa},
\mathbf{a}_{\mathcal{N}_{i}^\kappa},
\widetilde{\mathbf{a}}_{\mathcal{N}_{-i}^\kappa}
)| \\
\le c\rho^{\kappa+1} .
\end{split}
\end{equation}

This means the effect of different agents' observations and actions on any individual Q-function exponentially decreases as the topology distance increases, as illustrated in Fig. \ref{Exponential Decay Property}. Fortunately, this property holds for our system, as evidenced by the following proposition in our prior work \cite{quScalableReinforcementLearning2021}, which shows exponential decay property holds generally for networked MDPs.
\newtheorem{theorem}{Theorem}
\label{Thm: exponential decay property holds}

\begin{proposition}[Lemma 2 in \cite{quScalableReinforcementLearning2021}]
    Assume $\forall i$, $r_i$ is upper bounded by $\overline{r}$. Then the $(\frac{\overline{r}}{1-\gamma},\gamma)$-exponential decay property holds.
\end{proposition}

With this property, we further propose the SNA framework where the training of critics and actors leverages the network structure to achieve scalability.

\subsubsection{Critics that take only observations and actions of $\kappa$-hop neighbors}

Using the exponential decay property, we define a class of truncated Q-functions, which depend only on the states and actions of the local and $\kappa$-hop neighbors
\begin{equation}
\label{Truncated Q-function}
\begin{split}
\hat{Q}_{i}^{\pi}(\mathbf{o}_{\mathcal{N}_i^\kappa},\mathbf{a}_{\mathcal{N}_i^\kappa}) =
\mathbb{E}_{
\mathbf{o}_{\mathcal{N}_{-i}^\kappa},
\mathbf{a}_{\mathcal{N}_{-i}^\kappa}
}[Q_{i}^{\pi}
(\mathbf{o}_{\mathcal{N}_{i}^\kappa},
\mathbf{o}_{\mathcal{N}_{-i}^\kappa},
\mathbf{a}_{\mathcal{N}_{i}^\kappa},
\mathbf{a}_{\mathcal{N}_{-i}^\kappa}
)],
\end{split}
\end{equation}
where the expectation on $\mathbf{o}_{\mathcal{N}_{-i}^\kappa},
\mathbf{a}_{\mathcal{N}_{-i}^\kappa}$ can be over arbitrary distributions.

Due to the exponential decay property, such truncated individual Q-functions are still good approximations to the original individual Q-functions, because we have the following equation:
\begin{equation}
\label{Truncated Q-function error analysis}
\begin{array}{l}
|Q_i^\pi ({{\bf{o}}_{{\cal N}_i^\kappa }},{{\bf{o}}_{{\cal N}_{ - i}^\kappa }},{{\bf{a}}_{{\cal N}_i^\kappa }},{{\bf{a}}_{{\cal N}_{ - i}^\kappa }}) - \hat Q_i^\pi ({{\bf{o}}_{{\cal N}_i^\kappa }},{{\bf{a}}_{{\cal N}_i^\kappa }})|\\[1.5ex]
 \le {\mathbb{E}_{{{{\bf{\tilde o}}}_{{\cal N}_{ - i}^\kappa }},{{{\bf{\tilde a}}}_{{\cal N}_{ - i}^\kappa }}}}[|Q_i^\pi ({{\bf{o}}_{{\cal N}_i^\kappa }},{{\bf{o}}_{{\cal N}_{ - i}^\kappa }},{{\bf{a}}_{{\cal N}_i^\kappa }},{{\bf{a}}_{{\cal N}_{ - i}^\kappa }})\\[1.5ex]
 - Q_i^\pi ({{\bf{o}}_{{\cal N}_i^\kappa }},{{{\bf{\tilde o}}}_{{\cal N}_{ - i}^\kappa }},{{\bf{a}}_{{\cal N}_i^\kappa }},{{{\bf{\tilde a}}}_{{\cal N}_{ - i}^\kappa }})|]\\[1.5ex]
 \le c\rho^{\kappa+1}.
\end{array}
\end{equation}
Therefore, when using neural networks to approximate individual Q-functions, we can restrict their input variables to only the observations and actions of $\kappa$-hop neighbors. The above equation \eqref{Truncated Q-function error analysis} shows such input variables are adequately informative to approximate the true individual Q-functions. Importantly, such truncated Q-functions do not have their input dimension increase with the number of agents. This resolves the issue in the CTDE framework where critics' Q-functions need to input global information to estimate the global expected rewards in order to evaluate the corresponding actors' actions.

\subsubsection{Actors that are guided by critics of its $\kappa$-hop neighbors}
As mentioned above, each agent's critic uses truncated individual Q-function to estimate the expected individual rewards. As a result, if each actor is guided only by its own critic, the trained policy may overlook the global welfare. Therefore, to enhance cooperation among agents and improve global welfare, each local actor, in the proposed SNA framework, is guided by the critics of its $\kappa$-hop neighbors
\begin{equation}
\label{Eq: k-hop neighbors' critics}
    \sum_{j\in \mathcal{N}_i^\kappa}\hat{Q}_{j}^\pi(\mathbf{o}_{\mathcal{N}^{\kappa}_j,t},\mathbf{a}_{\mathcal{N}^{\kappa}_j,t}).
\end{equation}
In this way, agent $i$'s policy can increase the shared rewards of itself and its $\kappa$-hop neighbors.

Compared with the conventional CTDE framework shown in Fig. \ref{CTDE framework} (a), the proposed framework (with $\kappa$ equals 1), as shown in Fig. \ref{CTDE framework} (b), exploits the network structure and trains all the actors in a distributed way to achieve better scalability. 

\subsection{Multi-agent SAC algorithm with the SNA framework}

\begin{algorithm}[t]
\label{Alg1}
\caption{Multi-agent SAC with SNA framework}
\begin{algorithmic}
\State Initialize experience buffer $\mathcal{D}$, and parameters of the actors and critics.
\For{each episode}
    \For{each step $t$}
        \For{each agent $i$}{ in parallel}
            \State Choose $\mathbf{a}_{i,t}$ based on $\pi_{\phi_i}(\cdot|\mathbf{o}_{i,t})$;
            \State Take actions;
            \State Get reward $r_{i,t}$ and next states $\mathbf{s}_{i,t+1}$ from 
            \State the environment;
        \EndFor
        \State $\mathcal{D} \gets \mathcal{D}\cup \{\mathbf{s}_{t},\mathbf{a}_{t},r_{t},\mathbf{s}_{t+1}\}$
        \For{each agent $i$}{ in parallel}
            \State Sample a batch from experience buffer $\mathcal{D}$;
            \State Update individual Q-function: 
            \State $\theta_i \gets \theta_i+\sigma_i \nabla_{\theta_i}J_{Q_i}(\theta_i)$;
            \State Update individual policy: 
            \State $\phi_i \gets \phi_i + \sigma_i \nabla_{\phi_i}J_{\pi_i}(\phi_i)$;
            \State Softly update target networks' parameters;
        \EndFor
    \EndFor
\EndFor
\end{algorithmic}
\end{algorithm}

In this section, the MASAC algorithm \cite{haarnojaSoftActorCriticOffPolicy2018} is integrated with the proposed SNA framework to achieve scalability. Below, we explain how to train both the actor and critic.

In terms of the actors, the actors take a random action using  
\begin{equation}
\label{Policy Parameterization}
\mathbf{a}_{i,t}=\mathbf{f}_{\phi_i}(\mathbf{o}_{i,t};\mathbf{\epsilon}_t)
\end{equation}
where $\mathbf{\epsilon}_t$ is a random variable sampled from some fixed distribution, such as a spherical Gaussian. To be more specific, $\mathbf{f}_{\phi_i}$ is defined as
\begin{equation}
\label{Specific policy parameterization}
\begin{split}
    \mathbf{f}_{\phi_i}(\mathbf{o}_{i,t};\mathbf{\epsilon}_t) = \text{tanh}(\mu_{\phi_i}(\mathbf{o}_{i,t})+\sigma_{\phi_i}(\mathbf{o}_{i,t}) \odot \mathbf{\epsilon}_t), \\
    \mathbf{\epsilon}_t \sim \mathcal{N}(0,I).
\end{split}
\end{equation}
With $\mathbf{f}_{\phi_i}$ given, the policy $\pi_{\phi_i}(\cdot|\mathbf{o}_{i,t})$ is the probability distribution defined implicitly by $\mathbf{f}_{\phi_i}$. 

To train the critics, the soft truncated Q-function parameters are updated to minimize the following objective function
\begin{equation}
\label{Eq: SAC_critic}
\begin{split}   
J_{Q_i}(\theta_i)=\mathbb{E}_{\mathbf{x} \sim \mathcal{D}}[\frac{1}{2}(\hat{Q}_{\theta_i}(\mathbf{o}_{\mathcal{N}_{i}^\kappa,t},\mathbf{a}_{\mathcal{N}_{i}^\kappa,t})
-\hat{y})^2],
\end{split}
\end{equation}
with 
\begin{equation}
\label{Eq: SAC_critic_target}
\begin{split}
    \hat{y} = r_{i,t}+\gamma[\hat{Q}_{\overline{\theta}_i}(\mathbf{o}_{\mathcal{N}_{i}^\kappa,t+1},\mathbf{a'}_{\mathcal{N}_{i}^\kappa,t+1})\\
    -\alpha \log(\pi_{\phi_i}(\mathbf{a'}_{i,t+1}|\mathbf{o}_{i,t+1}))],
\end{split}
\end{equation}
where $\mathbf{x}$ denotes a tuple $(\mathbf{o}_{\mathcal{N}_{i}^\kappa,t},\mathbf{a}_{\mathcal{N}_{i}^\kappa,t},r_{i,t},\mathbf{o}_{\mathcal{N}_{i}^\kappa,t+1})$ sampled from a replay buffer $\mathcal{D}$; $\alpha$ is the temperature parameter; $\overline{\theta}_i$ is the parameter of the i-agent's target network which is softly updated using $\overline{\theta}_i\leftarrow \tau \theta_i + (1-\tau) \overline{\theta}_i$; $\mathbf{a'}_{\mathcal{N}_{i}^\kappa,t+1}$ is the sampled actions of agents in $\mathcal{N}_i^\kappa$ based on observations $\mathbf{o}_{\mathcal{N}_{i}^\kappa,t+1}$. 

Based on (\ref{Eq: SAC_critic_target}), it can be observed that the truncated Q-function merely estimates the expected individual rewards and the entropy of the policy adopted by the corresponding actor. Moreover, its input variables consist solely of the observations and actions of $\kappa$-hop neighbors.

To train the actors, the policy parameters are updated to minimize the following equations: 
\begin{equation}
\label{Eq: SNA actor objective}
\begin{split}
    J_{\pi_i}(\phi_i)=\mathbb{E}_{(\mathbf{o}_{t})\sim \mathcal{D}}[\alpha \log(\pi_{\phi_i}(\mathbf{a'}_{i,t}|\mathbf{o}_{i,t})-\\
    \frac{1}{N}\sum_{j\in \mathcal{N}_i^\kappa}\hat{Q}_{\theta_j}(\mathbf{o}_{\mathcal{N}^{\kappa}_j,t},\mathbf{a'}_{\mathcal{N}^{\kappa}_j,t})],
\end{split}
\end{equation}
where $\mathbf{a'}_{\mathcal{N}_{j}^\kappa,t+1}$ is the sampled actions of agents in $\mathcal{N}_j^\kappa$ based on observations $\mathbf{o}_{\mathcal{N}_{j}^\kappa,t+1}$. 

As a result, each actor is guided by the critics of its $\kappa$-hop neighbors, so that it can be trained to contribute to the welfare of itself and its $\kappa$-hop neighbors. 

Finally, the SAC algorithm with SNA framework is shown in Algorithm 1. Compared with the original SAC algorithm, this modified algorithm can realize better scalability, which is further discussed in the next subsection. 

\subsection{Rationale}
For the training of critics, according to (\ref{Truncated Q-function error analysis}), the truncated individual Q-functions provide a good estimation of the original Q-functions. 

For the training of actors, however, further explanation is needed to understand whether being guided only by the critics of the neighbors is sufficient to maximize the global welfare. To answer this question, we first assume that the real Q-functions are known, so that the actor objective in the CTDE framework can be rewritten as KL-divergence:
\begin{equation}
\label{Eq: Original KL}
\mathbb{E}_{(\mathbf{o}_{t})\sim \mathcal{D}}[\text{D}_{\text{KL}}(\pi_{\phi_i}(\cdot | \mathbf{o}_{i,t})||\frac{\text{exp}(\frac{1}{N}\sum_{j}Q_{j}^{\pi}(\mathbf{o}_t,\cdot,\mathbf{a'}_{-i,t}))}{Z(\frac{1}{N}\sum_{j}Q_{j}^{\pi}(\mathbf{o}_t,\cdot,\mathbf{a'}_{-i,t}))})]
\end{equation}
where the temperature parameter is omitted explicitly, as it can always be subsumed into the reward by scaling it by $\alpha^{-1}$, $\mathbf{a'}_{-i,t}$ are actions sampled from the current policies, and $Z(\sum_{j}Q_{j}^{\pi}(\mathbf{o}_t,\cdot,\mathbf{a}_{-i,t}))$ is the partition function that normalizes the distribution and is show as
\begin{equation}
\label{Eq: partitionfunction}
\begin{split}
&Z(\frac{1}{N}\sum_{j}Q_{j}^{\pi}(\mathbf{o}_t,\cdot,\mathbf{a'}_{-i,t}))\\
=&\int_{\mathcal{A}_i}\text{exp}(\frac{1}{N}\sum_{j}Q_{j}^{\pi}(\mathbf{o}_t,\mathbf{a}_{i,t},\mathbf{a'}_{-i,t}))d\mathbf{a}_{i,t}.
\end{split}
\end{equation}

As in the SNA framework where truncated Q-functions are used, the actor objective (\ref{Eq: SNA actor objective}) can be rewritten as the following KL-divergence:
\begin{equation}
\label{Eq: SNA KL}
\begin{split}
\mathbb{E}_{(\mathbf{o}_{t})\sim \mathcal{D}}[&\text{D}_{\text{KL}}
(
\pi_{\phi_i}(\cdot | \mathbf{o}_{i,t}) 
||\\
&\frac{
    \text{exp}(\frac{1}{N}\sum_{j\in \mathcal{N}_i^{\kappa}}\hat{Q}_{j}^{\pi}(\mathbf{o}_{\mathcal{N}_j,t},\cdot,\mathbf{a}_{\mathcal{N}_j/i,t}))
    }
    {Z(\frac{1}{N}\sum_{j\in \mathcal{N}_i^{\kappa}}\hat{Q}_{j}^{\pi}(\mathbf{o}_{\mathcal{N}_j,t},\cdot,\mathbf{a}_{\mathcal{N}_j/i,t}))
    }
)].
\end{split}
\end{equation}

Our goal is then to compare (\ref{Eq: Original KL}) and (\ref{Eq: SNA KL}).
If the objective in (\ref{Eq: Original KL}) is close to that in (\ref{Eq: SNA KL}), it implies that the actor objective utilizing solely neighboring critics within the proposed SNA framework aligns with the actor objective as if the non-truncated Q-function across all agents were known and used. 
This assertion would substantiate the proposed actor training within the SNA framework. To analyze the difference between (\ref{Eq: Original KL}) and (\ref{Eq: SNA KL}), note that the $(c,\rho)$-exponential decay property holds. Therefore, the error between the truncated Q-functions and the Q-function can be controlled, as shown in equation (\ref{Truncated Q-function error analysis}). This eventually leads to the following proposition.

\begin{proposition}
\label{Prop: Rationale behind actors training}
When (\ref{Exponential Decay Property}) and (\ref{Truncated Q-function error analysis}) holds, the absolute difference between the expectations of KL-divergences (\ref{Eq: Original KL}) and (\ref{Eq: SNA KL}) is bounded by $2c\rho^{\kappa+1}$.
\end{proposition}
The proof can be found in the Appendix. This proposition gives us confidence that the actors trained with truncated Q-functions can perform similarly with those trained with original Q-functions. Therefore, even if these actors are only guided by its $\kappa$-hop neighbors, they can perform to optimize the global welfare.     

\section{Case Study}
In this section, we detail the setup of our training environments and present the training results. We carry out several case studies with the number of agents varying from 6 to 114. This allows us to compare the training performance of our algorithm with that of other algorithms, demonstrating its superior scalability. Additionally, we explore key features of the Scalable Network-Aware (SNA) framework, including its compatibility with other actor-critic algorithms and the effect of the unique hyper-parameter $\kappa$ in our model.

\subsection{Experimental Settings}
\begin{figure*}[!t]
\centering
\includegraphics[width=7.25in]{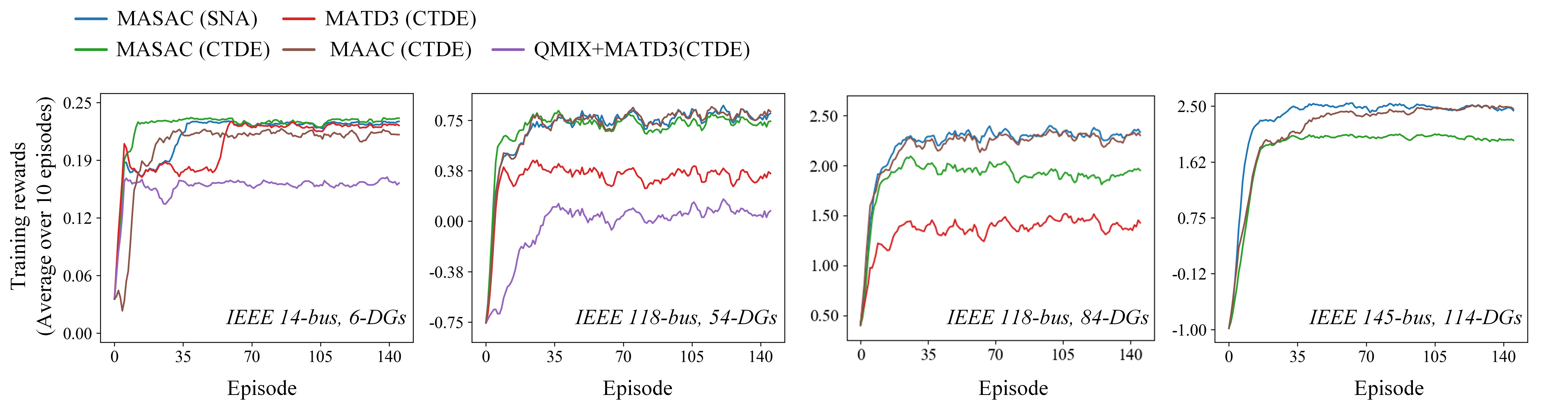}
   \caption{MARL algorithms training curves comparison.}
\label{fig: Training Results Comparison}
\end{figure*}

\begin{figure}[!t]
\centering
\includegraphics[width=3.25in]{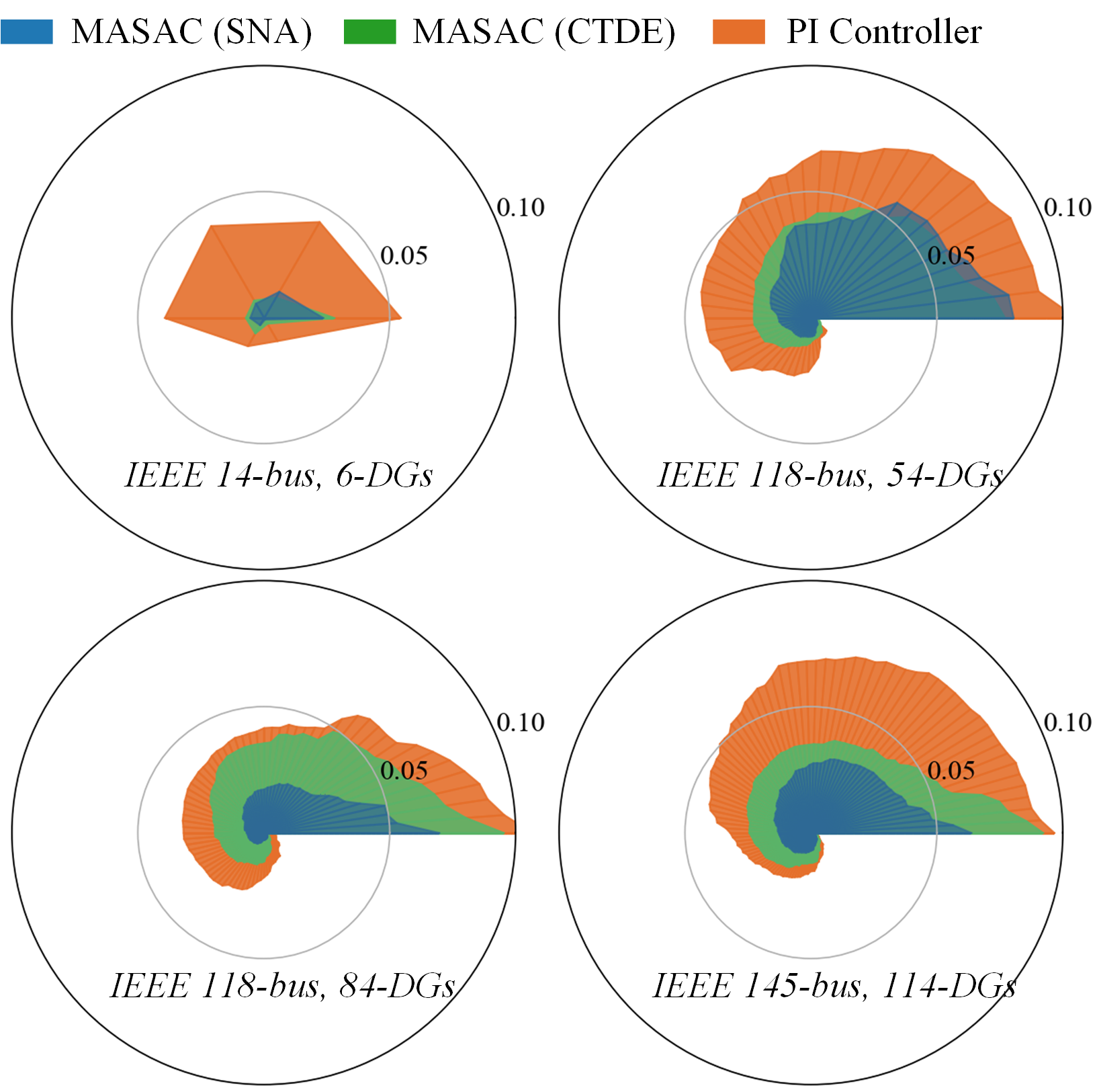}
   \caption{Voltage fluctuation comparison.}
\label{fig: Voltage fluctuation comparison}
\end{figure}

To construct the training environment, the transient models are adopted to simulate the dynamic behaviors. 
The detailed transient models of the DGs can be found in \cite{bidramDistributedCooperativeSecondary2013a}, while the models of the transmission lines and loads are illustrated in Fig. \ref{DER_load_transmission_model}. As mentioned, there are two levels of control for the DGs. The primary control is modeled as a part of the training environment, and the secondary control is implemented as an agent that utilizes neural networks. The control period is set to be 0.025s. As for the load, the active and reactive power undergo a random change within $\pm 5\%$ at each moment, a Markov process that mimics the random changes in real-world loads and non-dispatchable DGs. Furthermore, this environment supports the input of the standard data format of MatPower \cite{zimmermanMATPOWERSteadyStateOperations2011}, and simulate the dynamic behaviors based on the corresponding topology given in MatPower data.

In the experiments, we modify the original standard cases of MatPower by changing the number of DGs to test the scalability of our algorithm, and all deep learning algorithms are implemented using the PyTorch deep learning framework. The MAAC algorithm utilizes multi-head attention neural networks, whereas the neural networks employed in the other deep reinforcement learning algorithms are simple Feed-Forward Neural Networks (FFNNs). For algorithms that use FFNNs, the neural network structures are all identical whose the number of hidden layers is 2 and the hidden size is 256. 

\subsection{Training Results Comparisons}

In the experiment, four cases are set up, where the number of DGs (e.g., agents) are 6, 54, 84, and 114 respectively. The scalability of the proposed SNA framework is showcased by comparing it with the MASAC algorithm under the CTDE framework, and the MATD3 algorithm under the CTDE framework. Furthermore, comparisons are also made with the MAAC algorithm and the QMIX+MATD3 algorithm. The MAAC algorithm employs attention neural networks to encode the high-dimensional input for the critic, enabling the neural network to focus on more relative and significant information. The QMIX algorithm, originally a class of Q-learning, is combined with MATD3 in \cite{wangRealTimeJointRegulations2023} with the aim of addressing the credit assignment problem and improving scalability. 

The training rewards for each algorithm under the four cases are presented in Fig. \ref{fig: Training Results Comparison}. Note that the results for the SNA framework shown in this figure are given for $\kappa=1$. 

As can be observed, compared to the MASAC algorithm under the CTDE framework, MASAC under the SNA framework yields similar final training rewards in cases with 6 and 54 DGs. However, the SNA framework achieves higher final training rewards in cases with 84 and 114 DGs. This indicates that as the number of DGs increases, the proposed SNA framework exhibits improved training performance. This improvement is attributed to the truncation of critics' input by leveraging the network structure, thereby retaining only the most influential input and achieving better scalability. The performance of the QMIX+MATD3 algorithm is suboptimal in all cases, possibly due to the QMIX algorithm's assumption that each agent's individual policy can improve global rewards by boosting individual rewards, an assumption that may not be easily satisfied in our test cases. Lastly, in examining the MAAC algorithm, its training results are comparable to those of the proposed SNA framework. 
However, it is worth noting that the attention network is considerably more complex than the Feed-Forward Neural Network (FFNN) we employed, resulting in a larger computational burden. For a fair comparison with MAAC, the SNA framework should also use the attention network (as opposed to the simple FFNN), and this is left as a future direction. 

To show the effectiveness and advantages of the proposed SNA framework in regulating the voltage, we compare the node voltage fluctuations under the control of MASAC under the SNA framework, MASAC under the CTDE framework, and traditional Proportional Integral (PI) control, as shown in Fig \ref{fig: Voltage fluctuation comparison}. The figures are plotted through the following process. Run the simulation for 100 time points, record the maximum per-unit fluctuation value reached by the node, and then sort and connect these fluctuation values in a clockwise direction from the largest to the smallest, thereby visually displaying the fluctuation range. Hence, a smaller area implies less voltage fluctuation under the respective control. As can be observed, the MASAC algorithm under both the SNA and CTDE frameworks perform much better than the traditional PI controllers. In cases with 84 and 114 DGs, the SNA framework significantly reduces voltage fluctuations and decreases the proportion of per-unit voltage fluctuations exceeding 0.05. 

\subsection{Influence of $\kappa$ Value}

\begin{figure}[!t]
\centering
\includegraphics[width=2.75in]{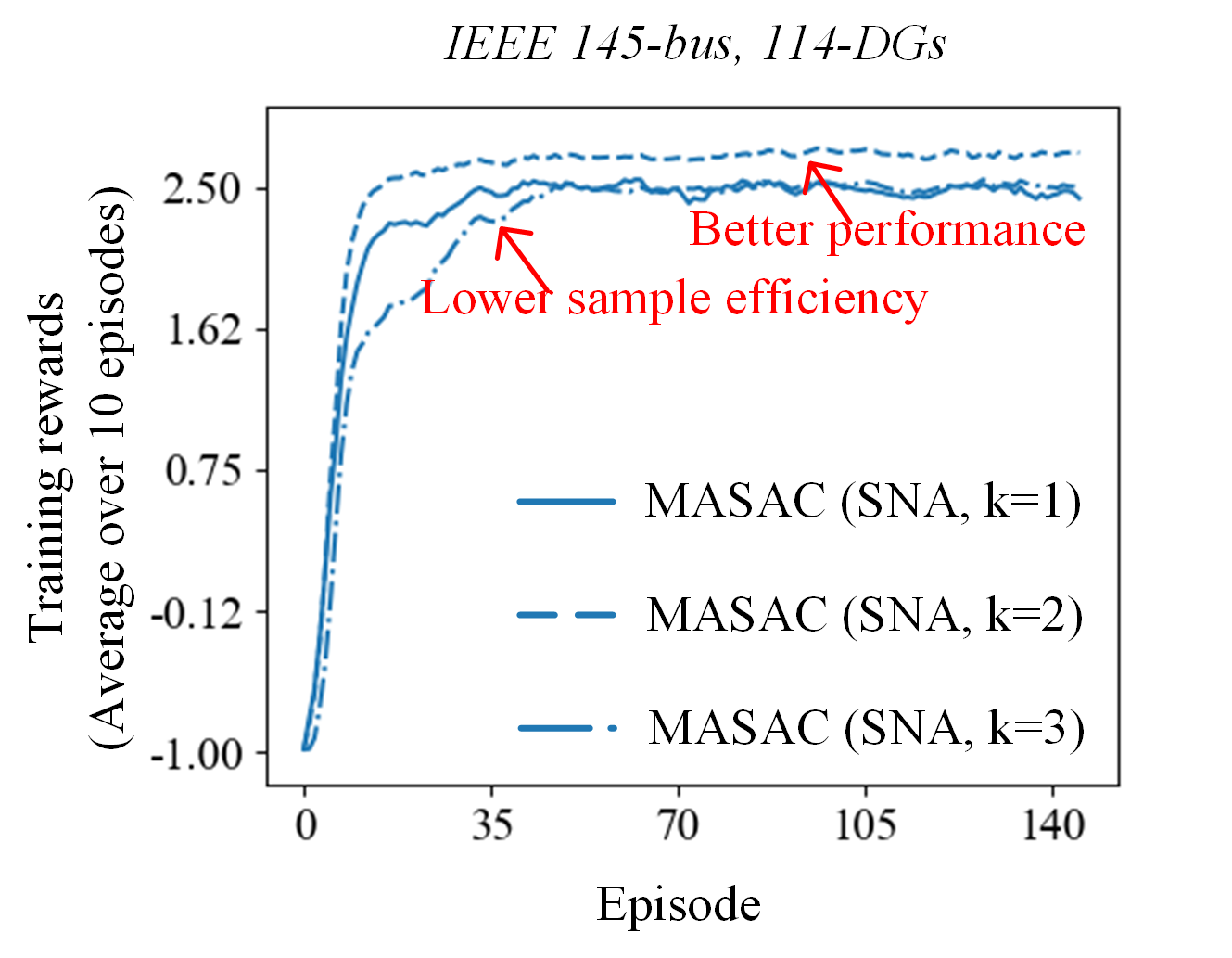}
   \caption{Training rewards under different $\kappa$ values.}
\label{fig: kappaInfluence}
\end{figure}

\begin{figure}[!t]
\centering
\includegraphics[width=2.75in]{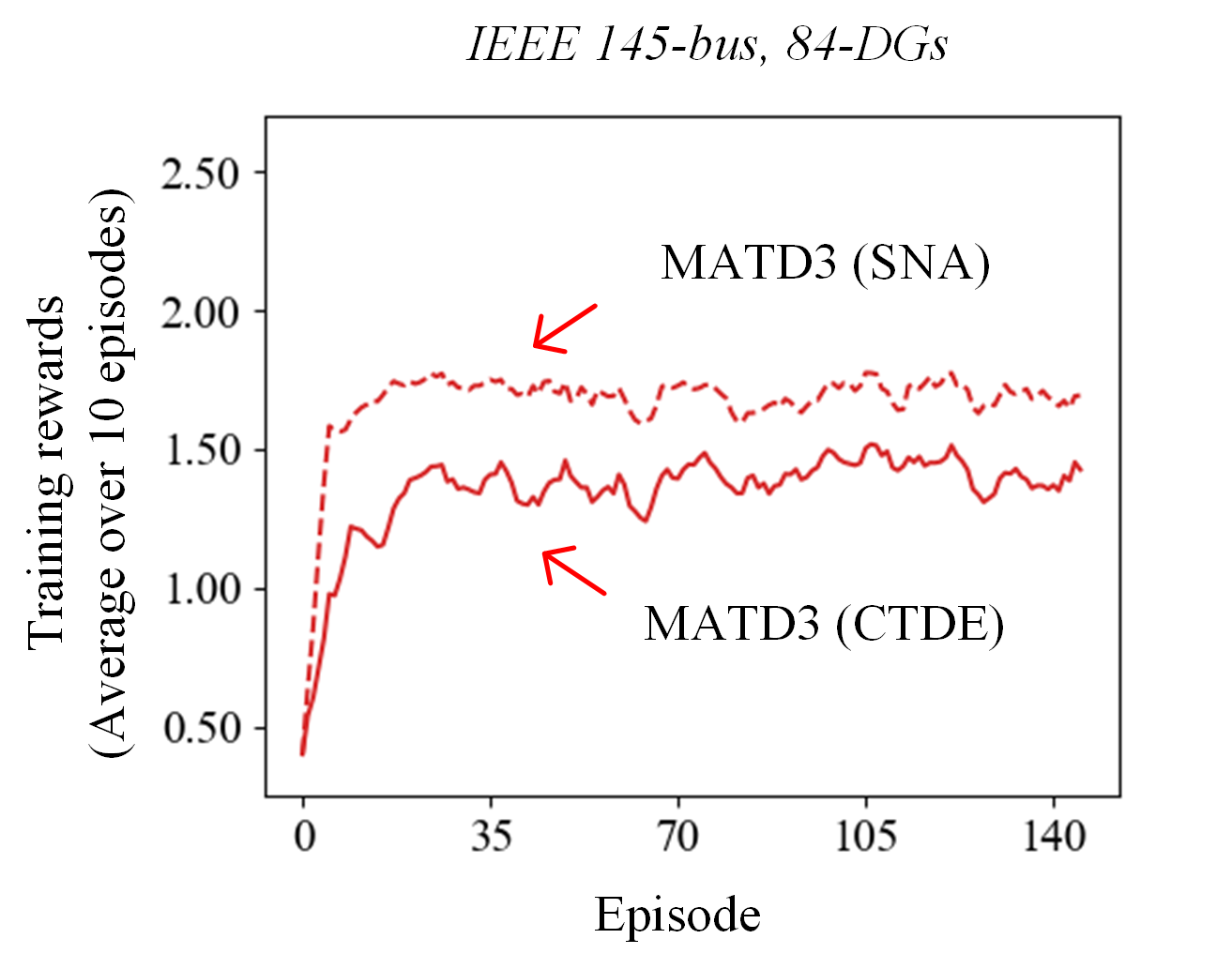}
   \caption{Compatibility with MATD3 algorithm.}
\label{fig: compatibility with MATD3}
\end{figure}

\begin{table}[!t]
\centering
\caption{Input Dimension Comparisons}
\begin{tabular}{cccc}
\hline\hline
\multirow{2}{*}{Case} & \multirow{2}{*}{MARL Framework} & \multicolumn{2}{c}{Critics input dimension}\\
                      &                & Maximum & Average \\
\hline
\multirow{2}{*}{Node14-6} & CTDE & 60 & 60\\
                            & SNA ($\kappa=1$) &  30 & 20\\
\multirow{2}{*}{Node118-54} & CTDE & 540 & 540\\
                            & SNA ($\kappa=1$) &  60 & 23\\
\multirow{2}{*}{Node118-84} & CTDE & 840 & 840\\
                            & SNA ($\kappa=1$) &  80 & 30\\
\multirow{4}{*}{Node145-114} & CTDE & 1140 & 1140\\
                            & SNA ($\kappa=1$) &  180 & 54\\
                            & SNA ($\kappa=2$) &  500 & 186\\
                            & SNA ($\kappa=3$) &  770 & 381\\
\hline\hline
\end{tabular}
\label{table:1}
\end{table}

The hyper-parameter $\kappa$ determines how many neighbors' observations and action are taken by the critics to estimate expected individual rewards and how many neighbor's critics are used to guide the actors to improve the welfare of these neighbors. According to (\ref{Truncated Q-function error analysis}), as $\kappa$ increases, the truncation error of the truncated Q-function becomes smaller, thereby enabling a more accurate estimation of individual rewards. However, an increase in $\kappa$ leads to an increase in input dimension as shown in Table \ref{table:1}, degrading training performance due to the curse of dimensionality.

In the 114-DG case, the value of $\kappa$ is varied to show its influence on training performance. The training results are shown in Fig. \ref{fig: kappaInfluence}. When $\kappa$ equals 2, the SNA framework achieves greater final training rewards. However, when $\kappa$ equals 3, the training performance of the SNA framework is compromised due to an increase in input dimension, leading to lower sample efficiency.

\subsection{Compatibility with Other Actor-Critic Algorithms}

As discussed in the introduction, the proposed SNA framework can seamlessly integrate with multi-agent actor-critic algorithms. So far, all the implementation of SNA has been based on MASAC. In this section, the proposed framework is combined with the MATD3 algorithm and trained in the 84-DG case. The training results are shown in Fig. \ref{fig: compatibility with MATD3}. It can be seen that the proposed framework has improved the final training rewards of the algorithm. This demonstrates the compatibility of the proposed framework with other actor-critic type algorithms.

\subsection{Comprehensive Evaluation}
The innovation and key to achieving excellent scalability of the proposed method lies in leveraging the network structure. Firstly, from the perspective of problem formulation, the voltage problem is formulated as a networked MDP. This modeling approach, distinct from traditional methods, retains information about the network structure. Secondly, by leveraging the exponential decay property in the networked MDP, the inputs to the critics are truncated, and the actors are trained under the guidance of its neighbors' critics based on the network structure. As a result, the critics can be effectively trained and the actors can be guided to optimize the global welfare even under the cases with large number of DGs. 

Another benefit associated with the SNA framework is its facilitation of distributed training. The hyper-parameter $\kappa$ in the SNA framework determines the amount of information that critics and actors need to acquire during training. In the studied cases, when $\kappa=1$, our proposed framework can achieve superior performance. This means only communications between neighboring DGs are required, greatly ensuring the practicality of the proposed method.

\section{Conclusion}
In this article, the SNA framework is proposed to achieve superior scalability for implementing decentralized inverter-based voltage control. This advantage is realized by utilizing the network structure to effectively train both critics and actors. Moreover, this framework facilitates distributed training as it solely relies on neighboring information. 

In the studied cases of this article, the SNA framework is combined with MASAC and MATD3 algorithms. The incorporation of the MASAC algorithm within the SNA framework facilitates effective voltage control in a system encompassing 114 DGs which is the state-of-the-art result. In future work, this method holds potential for application in other power system control scenarios, such as frequency control.

{\appendix[Proof of Proposition \ref{Prop: Rationale behind actors training}]
We introduce the following notations to simplify the expressions, and for simplicity, we just use $f$, $\overline{f}$, and $\hat{f}$ for notation when not causing misunderstandings:

\begin{equation}
\label{Eq: denotions}
\left\{
\begin{aligned}
&f(\mathbf{o}_{t},\mathbf{a}_{i,t},\mathbf{a}_{-i,t}) = \frac{1}{N}\sum_{j}Q_{j}^\pi(\mathbf{o}_{t},\mathbf{a}_{i,t},\mathbf{a}_{-i,t})\\
&\overline{f}(\mathbf{o}_{t},\mathbf{a}_{i,t},\mathbf{a}_{-i,t})=\frac{1}{N}\sum_{j\in \mathcal{N}_i^\kappa}Q_{j}^\pi(\mathbf{o}_{t},\mathbf{a}_{i,t},\mathbf{a}_{-i,t})\\
&\hat{f}(\mathbf{o}_{t},\mathbf{a}_{i,t},\mathbf{a}_{-i,t})=\frac{1}{N}\sum_{j\in \mathcal{N}_i^\kappa}\hat{Q}_{j}^\pi(\mathbf{o}_{\mathcal{N}_j^\kappa,t},\mathbf{a}_{i,t},\mathbf{a}_{\mathcal{N}_j^\kappa/i,t})
\end{aligned}
\right. .
\end{equation}

Firstly, we examine the differences among these three functions, for which the following notations are introduced:
\begin{equation}
\label{Eq: Appendix4}
\begin{split}
\Delta\overline{f}=f-\overline{f}
=\frac{1}{N}\sum_{j\notin \mathcal{N}_i^\kappa}Q_{j}^\pi(\mathbf{o}_{t},\mathbf{a}_{i,t},\mathbf{a}_{-i,t}),
\end{split}
\end{equation}
\begin{equation}
\label{Eq: Appendix4_1}
\begin{split}
\Delta\hat{f}=\overline{f}-\hat{f}
=&\frac{1}{N}\sum_{j\in \mathcal{N}_i^\kappa}Q_{j}^\pi(\mathbf{o}_{t},\mathbf{a}_{i,t},\mathbf{a}_{-i,t})-\\
&\hat{Q}_{j}^\pi(\mathbf{o}_{\mathcal{N}_j^\kappa,t},\mathbf{a}_{i,t},\mathbf{a}_{\mathcal{N}_j^\kappa/i,t}).
\end{split}
\end{equation}

According to the exponential decay property as shown in (\ref{Eq: Exponential decay property}), the following inequality can be derived:
\begin{equation}
\label{Eq: Appendix5}
\begin{split}
\max_{\mathbf{a}_{i,t}}\Delta\overline{f}-\min_{\mathbf{a}_{i,t}}\Delta\overline{f}\le (1-\frac{|\mathcal{N}_i^\kappa|}{N})c\rho^{\kappa+1}.
\end{split}
\end{equation}

According to the truncated errors of truncated individual Q-functions as shown in (\ref{Truncated Q-function error analysis}), the following inequality can be derived:
\begin{equation}
\label{Eq: Appendix7}
\begin{split}
\max_{\mathbf{a}_{i,t}}\Delta\hat{f}-\min_{\mathbf{a}_{i,t}}\Delta\hat{f}\le 2\frac{|\mathcal{N}_i^\kappa|}{N}c\rho^{\kappa+1}.
\end{split}
\end{equation}

Then, we examine the difference between (\ref{Eq: Original KL}) and (\ref{Eq: SNA KL}), according to the definition of KL divergence, the following relation can be obtained:
\begin{equation}
\label{Eq: Appendix9}
\begin{split}
&|\text{D}_{\text{KL}}(\pi_{\phi_i}(\cdot | \mathbf{o}_{i,t})||\frac{e^f}{Z(f)})-\text{D}_{\text{KL}}(\pi_{\phi_i}(\cdot | \mathbf{o}_{i,t})||\frac{e^{\hat{f}}}{Z(\hat{f})})|\\
=&|-\int_{\mathcal{A}_i}\pi_{\phi_i}(\mathbf{a}_{i,t}|\mathbf{o}_{i,t})\ln(\frac{e^fZ({\overline{f}})}{e^{\overline{f}}Z(f)})d\mathbf{a}_{i,t}-\\
&\int_{\mathcal{A}_i}\pi_{\phi_i}(\mathbf{a}_{i,t}|\mathbf{o}_{i,t})\ln(\frac{e^{\overline{f}}Z({\hat{f}})}{e^{\hat{f}}Z({\overline{f}})})d\mathbf{a}_{i,t}|,
\end{split}
\end{equation}
where the first and second terms on the right-hand side of the above equations can be bounded, as shown below:
\begin{equation}
\label{Eq: Appendix6}
\begin{split}
&|\int_{\mathcal{A}_i}\pi_{\phi_i}(\mathbf{a}_{i,t}|\mathbf{o}_{i,t})\ln(\frac{e^fZ({\overline{f}})}{e^{\overline{f}}Z(f)})d\mathbf{a}_{i,t}|\\
=&|\int_{\mathcal{A}_i}\pi_{\phi_i}(\mathbf{a}_{i,t}|\mathbf{o}_{i,t})\ln(\frac{e^{\Delta\overline{f}}Z({\overline{f}})}{Z({\Delta\overline{f}}+{\overline{f}})})d\mathbf{a}_{i,t}|\\
\le&\int_{\mathcal{A}_i}\pi_{\phi_i}(\mathbf{a}_{i,t}|\mathbf{o}_{i,t})\ln(\frac{e^{\max_{\mathbf{a}_{i,t}}\Delta\overline{f}}Z({\overline{f}})}{e^{\min_{\mathbf{a}_{i,t}}\Delta\overline{f}}Z({{\overline{f}})})}d\mathbf{a}_{i,t}\\
\le&(1-\frac{|\mathcal{N}_i^\kappa|}{N})c\rho^{\kappa+1},
\end{split}
\end{equation}
where the first inequality is obtained from the monotonicity of the logarithmic and exponential functions, while the second inequality is derived from (\ref{Eq: Appendix5}). Similarly, the upper limit of the second term is bounded by:
\begin{equation}
\label{Eq: Appendix8}
\begin{split}
&|\int_{\mathcal{A}_i}\pi_{\phi_i}(\mathbf{a}_{i,t}|\mathbf{o}_{i,t})\ln(\frac{e^{\overline{f}}Z({\hat{f}})}{e^{\hat{f}}Z({\overline{f}})})d\mathbf{a}_{i,t}|\\
=&|\int_{\mathcal{A}_i}\pi_{\phi_i}(\mathbf{a}_{i,t}|\mathbf{o}_{i,t})\ln(\frac{e^{\Delta\hat{f}}Z({\hat{f}})}{Z({\Delta\hat{f}}+{\hat{f}})})d\mathbf{a}_{i,t}|\\
\le&\int_{\mathcal{A}_i}\pi_{\phi_i}(\mathbf{a}_{i,t}|\mathbf{o}_{i,t})\ln(\frac{e^{\max_{\mathbf{a}_{i,t}}\Delta\hat{f}}Z({\hat{f}})}{e^{\min_{\mathbf{a}_{i,t}}\Delta\hat{f}}Z({{\hat{f}})})}d\mathbf{a}_{i,t}\\
\le&2\frac{|\mathcal{N}_i^\kappa|}{N}c\rho^{\kappa+1}.
\end{split}
\end{equation}

Combining (\ref{Eq: Appendix9}), (\ref{Eq: Appendix6}) and (\ref{Eq: Appendix8}), the following relationship can be obtained:
\begin{equation}
\label{Eq: Difference of KL divergence}
\begin{split}
&|\text{D}_{\text{KL}}(\pi_{\phi_i}(\cdot | \mathbf{o}_{i,t})||\frac{e^f}{Z(f)})-\text{D}_{\text{KL}}(\pi_{\phi_i}(\cdot | \mathbf{o}_{i,t})||\frac{e^{\hat{f}}}{Z(\hat{f})})|\\
\le & (1+\frac{|\mathcal{N}_i^\kappa|}{N})c\rho^{\kappa+1}\\
\le & 2c\rho^{\kappa+1},
\end{split}
\end{equation}
which completes the proof with the Jensen's inequality.

}

\bibliographystyle{IEEEtran}
\bibliography{ref}
 

\vfill

\end{document}